\documentclass{article}
\usepackage{graphicx}
\usepackage{amsmath}
\usepackage{amsfonts}
\usepackage{amssymb}
\newtheorem{theorem}{Theorem}

\newtheorem{corollary}[theorem]{Corollary}

\newtheorem{definition}[theorem]{Definition}

\newtheorem{lemma}[theorem]{Lemma}

\newtheorem{proposition}[theorem]{Proposition}

\newenvironment{proof}[1][Proof]{\textbf{#1.} }{\ \rule{0.5em}{0.5em}}

\begin{document}

\title{Singular polynomials for the symmetric group and Krawtchouk polynomials}
\author{Charles F. Dunkl\thanks{During the preparation of this paper the author was
partially supported by NSF grant DMS 0100539.}}
\date{16 October 2003}
\maketitle
\begin{abstract}
A singular polynomial is one which is annihilated by all Dunkl operators for a
certain parameter value. These polynomials were first studied by Dunkl, de Jeu
and Opdam, (\textit{Trans. Amer. Math. Soc.} 346 (1994), 237-256). This paper
constructs a family of such polynomials associated to the irreducible
representation $\left(  N-2,1,1\right)  $ of the symmetric group $S_{N}$ for
odd $N$ and parameter values $-\frac{1}{2},-\frac{3}{2},-\frac{5}{2},\ldots$.
The method depends on the use of Krawtchouk polynomials to carry out a change
of variables in a generating function involved in the construction of
nonsymmetric Jack polynomials labeled by $\left(  m,n,0,\ldots,0\right)
,m\geq n$.
\end{abstract}

\section{Introduction}

We will study polynomials on $\mathbb{R}^{N}$ with certain properties relating
to the action of the symmetric group $S_{N}$ acting\ as a finite reflection
group (of type $A_{N-1}$). Let $\mathbb{N}$ denote $\left\{  1,2,3,\ldots
\right\}  $ and $\mathbb{N}_{0}=\mathbb{N\cup}\left\{  0\right\}  $; for
$\alpha\in\mathbb{N}_{0}^{N}$ let $\left|  \alpha\right|  =\sum_{i=1}%
^{N}\alpha_{i}$ and define the monomial $x^{\alpha}$ to be $\prod_{i=1}%
^{N}x_{i}^{\alpha_{i}}$; its degree is $\left|  \alpha\right|  .$ Consider
elements of $S_{N}$ as functions on $\{1,2,\ldots,N\}$ then for $x\in
\mathbb{R}^{N}$ and $w\in S_{N}$ let $\left(  xw\right)  _{i}=x_{w\left(
i\right)  }$ for $1\leq i\leq N$; and extend this action to polynomials by
$wf\left(  x\right)  =f\left(  xw\right)  $. This has the effect that
monomials transform to monomials, $w\left(  x^{\alpha}\right)  =x^{w\alpha}$
where $\left(  w\alpha\right)  _{i}=\alpha_{w^{-1}\left(  i\right)  }$ for
$\alpha\in\mathbb{N}_{0}^{N}$. (Consider $x$ as a row vector, $\alpha$ as a
column vector, and $w$ as a permutation matrix, with $1$'s at the $\left(
w\left(  j\right)  ,j\right)  $ entries.) The reflections in $S_{N}$ are the
transpositions, denoted by $\left(  i,j\right)  $ for $i\neq j$, interchanging
$x_{i}$ and $x_{j}$.

In \cite{D1} the author constructed for each finite reflection group a
parametrized commutative algebra of differential-difference operators; for the
symmetric group there is one parameter $\kappa\in\mathbb{C}$ and the
definition is as follows:

\begin{definition}
For any polynomial $f$ on $\mathbb{R}^{N}$ and $1\leq i\leq N$ let
\[
\mathcal{D}_{i}f\left(  x\right)  =\frac{\partial}{\partial x_{i}}f\left(
x\right)  +\kappa\sum_{j\neq i}\frac{f\left(  x\right)  -\left(  i,j\right)
f\left(  x\right)  }{x_{i}-x_{j}}.
\]
\end{definition}

It was shown in \cite{D1} that $\mathcal{D}_{i}\mathcal{D}_{j}=\mathcal{D}%
_{j}\mathcal{D}_{i}$ for $1\leq i,j\leq N$ and each $\mathcal{D}_{i}$ maps
homogeneous polynomials to homogeneous polynomials. A specific parameter value
$\kappa$ is said to be a \textit{singular value} (associated with $S_{N}$) if
there exists a nonzero polynomial $p$ such that $\mathcal{D}_{i}p=0$ for
$1\leq i\leq N$; and $p$ is called a \textit{singular polynomial}. It was
shown by Dunkl, de Jeu and Opdam \cite[p.248]{DJO} that the singular values
for $S_{N}$ are the numbers $-\frac{j}{n}$ where $n=2,\ldots,N,\,j\in
\mathbb{N}$ and $n\nmid j$ ($n$ does not divide $j$). In this paper we
construct singular polynomials for the values -$\frac{j}{N-1}$ ($j\in
\mathbb{N}$ and $N-1\nmid j$) with a new result for the case of $N$ being odd
and -$\frac{j}{N-1}=-l-\frac{1}{2},l\in\mathbb{N}_{0}$. There are conjectures
in \cite[p.255]{DJO} regarding some general properties of the singular
polynomials for $S_{N}$ but these are not as yet established. Here is an easy
example of singular polynomials: let $a_{N}\left(  x\right)  =\prod
\limits_{1\leq i<j\leq N}\left(  x_{i}-x_{j}\right)  $, the alternating
polynomial. Then $\left(  i,j\right)  a_{N}\left(  x\right)  =-a_{N}\left(
x\right)  $ for any transposition ($i\neq j$); further $\frac{\partial
}{\partial x_{i}}a_{N}\left(  x\right)  =a_{N}\left(  x\right)  \sum_{j\neq
i}\frac{1}{x_{i}-x_{j}}$. Thus for any $l\in\mathbb{N}_{0}$ we have
$\mathcal{D}_{i}\left(  a_{N}\left(  x\right)  ^{2l+1}\right)  =\left(
2l+1+2\kappa\right)  \,a_{N}\left(  x\right)  ^{2l+1}\sum_{j\neq i}\frac
{1}{x_{i}-x_{j}}$ (for each $i$), which shows that $a_{N}\left(  x\right)
^{2l+1}$ is singular for $\kappa=-l-\frac{1}{2}$. Irreducible representations
of $S_{N}$ are labeled by partitions of $N$ (see, for example, Macdonald
\cite[p.114]{M}); the polynomial $a_{N}\left(  x\right)  ^{2l+1}$ is
associated with the representation $\left(  1,1,\ldots,1\right)  $ (more
precisely, the span $\mathbb{R}a_{N}\left(  x\right)  ^{2l+1}$ is an $S_{N}%
$-module of isotype $\left(  1,1,\ldots,1\right)  $).

Our construction is in terms of nonsymmetric Jack polynomials which are
defined to be the simultaneous eigenfunctions of the pairwise commuting
operators $\mathcal{D}_{i}x_{i}-\kappa\sum_{j<i}\left(  i,j\right)  $,$\,1\leq
i\leq N$ (details about these may be found in the book by Dunkl and Xu
\cite[Ch.8]{DX}). When $\kappa>0$ these operators are self-adjoint with
respect to the inner product on polynomials defined by
\[
\left\langle f,g\right\rangle _{\mathbb{T}}=c_{\kappa}\int_{\mathbb{T}^{N}%
}f\left(  x\right)  \,\overline{g\left(  x\right)  }\prod_{1\leq i<j\leq
N}\left|  x_{i}-x_{j}\right|  ^{2\kappa}dm\left(  x\right)  ,
\]
where $\mathbb{T}^{N}$ is the $N$-fold complex torus $\left\{  z\in
\mathbb{C}:\left|  z\right|  =1\right\}  ^{N},\,x_{j}=e^{\mathrm{i}\theta_{j}%
}$ for $-\pi<\theta_{j}\leq\pi$ (and $1\leq j\leq N$) and the standard measure
is $dm\left(  x\right)  =\prod_{j=1}^{N}d\theta_{j}$. The constant $c_{\kappa
}$ is chosen so that $\left\langle 1,1\right\rangle _{\mathbb{T}}=1$ (computed
by means of the Macdonald-Mehta-Selberg integral). The nonsymmetric Jack
polynomials are labeled by $\mathbb{N}_{0}^{N}$; in this paper only the labels
$\left(  m,n,0,\ldots,0\right)  $ will occur.

\section{The Basic Polynomials}

These are the relevant results from Dunkl \cite{D2}. All the polynomials
considered here have coefficients in $\mathbb{Q}\left(  \kappa\right)  $
(rational functions of $\kappa$ with rational coefficients). The polynomials
$p_{mn}\left(  x\right)  $ are defined by the generating function
\begin{equation}
\sum_{m,n=0}^{\infty}p_{mn}\left(  x\right)  s^{m}t^{n}=\left(  1-sx_{1}%
\right)  ^{-1}\left(  1-tx_{2}\right)  ^{-1}\prod_{i=1}^{N}\left(  \left(
1-sx_{i}\right)  \left(  1-tx_{i}\right)  \right)  ^{-\kappa} \label{genfun}%
\end{equation}
(for\ convergence require $\left|  s\right|  ,\left|  t\right|  <\left(
\max_{i}\left|  x_{i}\right|  \right)  ^{-1}$). Note that $p_{mn}\left(
x\right)  =x_{1}^{m}x_{2}^{n}$ when $\kappa=0$. In \cite{D2} it was shown that
$\mathcal{D}_{i}p_{mn}=0$ for all $i>2$, $\mathcal{D}_{2}p_{m0}=0,$
$\mathcal{D}_{1}p_{0n}=0$ and
\begin{equation}
\mathcal{D}_{1}p_{mn}=\left(  N\kappa+m\right)  p_{m-1,n}+\kappa\sum
_{i=0}^{n-1}\left(  p_{m+n-1-i,i}-p_{i,m+n-1-i}\right)  \label{D1p}%
\end{equation}
for $m>n\geq1$, and
\begin{align}
\mathcal{D}_{1}p_{mn}  &  =\left(  \left(  N-1\right)  \kappa+m\right)
p_{m-1,n}+\kappa p_{n,m-1}\label{D2p}\\
&  +\kappa\sum_{i=0}^{m-2}\left(  p_{m+n-1-i,i}-p_{i,m+n-1-i}\right) \nonumber
\end{align}
for $n\geq m\geq1$ (the second sum is omitted if $m=1$)$.$ The expression for
$\mathcal{D}_{2}p_{mn}$ can be deduced by interchanging the labels on $p$. The
Pochhammer symbol is $\left(  a\right)  _{n}=\prod_{i=1}^{n}\left(
a+i-1\right)  $ for any $a\in\mathbb{Q}\left(  \kappa\right)  $ and
$n\in\mathbb{N}_{0}$.

\begin{definition}
\label{wdef}For $m\geq n$ let
\begin{align*}
\omega_{mn}  &  =p_{mn}+\sum_{j=1}^{n}\frac{\left(  -\kappa\right)
_{j}\left(  m-n+1\right)  _{j-1}}{\left(  \kappa+m-n+1\right)  _{j}%
\,j!}\left(  \left(  m-n+j\right)  p_{m+j,n-j}+jp_{n-j,m+j}\right)  ,\\
\omega_{nm}  &  =\left(  1,2\right)  \omega_{mn}.
\end{align*}
\end{definition}

Observe that $\left(  1,2\right)  p_{mn}=p_{nm}$, and the coefficients of
$\omega_{mn}$ in the formula are independent of $N$. Also there is the
symmetry property $\left(  i,j\right)  \omega_{mn}=\omega_{mn}$ for $2<i<j\leq
N$. It can be shown that $\mathcal{D}_{1}x_{1}\omega_{mn}=\left(  \left(
N-1\right)  \kappa+m+1\right)  \omega_{mn}$ and $\left(  \mathcal{D}_{2}%
x_{2}-\kappa\left(  1,2\right)  \right)  \omega_{mn}=\left(  \left(
N-2\right)  \kappa+n+1\right)  \omega_{mn}$ for $m\geq n$. By means of the
product rule:
\begin{equation}
\mathcal{D}_{i}\left(  fg\right)  =f\mathcal{D}_{i}g+\frac{\partial
f}{\partial x_{i}}g+\kappa\sum_{j\neq i}\frac{f-\left(  i,j\right)  f}%
{x_{i}-x_{j}}\left(  i,j\right)  g, \label{prodrule}%
\end{equation}
we can show $\mathcal{D}_{i}x_{i}\omega_{mn}=\left(  1+\left(  N-3\right)
\kappa\right)  \omega_{mn}+\kappa\left(  \left(  1,i\right)  +\left(
2,i\right)  \right)  \omega_{mn}$ and\newline $\left(  \mathcal{D}_{i}%
x_{i}-\kappa\sum_{j<i}\left(  i,j\right)  \right)  \omega_{mn}=\left(  \left(
N-i\right)  \kappa+1\right)  \omega_{mn}$ for each $i>2$. Thus $\omega_{mn}$
is the nonsymmetric Jack polynomial labeled by $\left(  m,n,0,\ldots,0\right)
$). Note that $\omega_{mn}$ is defined whenever $\kappa\notin-\mathbb{N}$,
although $\omega_{mn}=0$ for certain values of $N,m,n$ and $\kappa.$ In fact
this is the key ingredient of our construction.

\begin{theorem}
\label{Dwmn}The following hold for all $\kappa\notin-\mathbb{N}$:

\begin{enumerate}
\item  for $m>n$, $\mathcal{D}_{1}\omega_{mn}=\left(  N\kappa+m\right)
\omega_{m-1,n}$\newline $+\dfrac{\left(  \left(  N-1\right)  \kappa+n\right)
\kappa}{\kappa+m-n}\left(  \omega_{n-1,m}-\dfrac{\kappa}{\kappa+m-n+1}%
\omega_{m,n-1}\right)  ;$

\item  for $m\geq n$, $\mathcal{D}_{2}\omega_{mn}=\left(  \left(  N-1\right)
\kappa+n\right)  \left(  \omega_{m,n-1}-\dfrac{\kappa}{\kappa+m-n+1}%
\omega_{n-1,m}\right)  ;$

\item  for $m=n$, $\mathcal{D}_{1}\omega_{nn}=\left(  \left(  N-1\right)
\kappa+n\right)  \left(  \omega_{n-1,n}-\dfrac{\kappa}{\kappa+1}\omega
_{n,n-1}\right)  .$
\end{enumerate}
\end{theorem}

\begin{proof}
It was shown in \cite[p.192]{D2} that both $\mathcal{D}_{1}\omega_{mn}$ and
$\mathcal{D}_{2}\omega_{mn}$ are in the span of $\left\{  \omega
_{m-1,n},\omega_{n,m-1},\omega_{m,n-1},\omega_{n-1,m}\right\}  .$ This implies
that only the coefficients of $p_{m-1,n},p_{n,m-1},p_{m,n-1},p_{n-1,m}$ need
to be calculated. Let $g_{1},g_{2},\ldots$ denote polynomials of the form
$\sum_{j=0}^{n-2}\left(  c_{j}p_{m+n-1-j,j}+c_{j}^{\prime}p_{j,m+n-1-j}%
\right)  $ with coefficients $c_{j},c_{j}^{\prime}\in\mathbb{Q}\left(
\kappa\right)  .$ Both formulae (\ref{D1p}), (\ref{D2p}) are used.

Suppose $m>n$, then
\begin{align*}
\mathcal{D}_{1}\omega_{mn}  &  =\left(  N\kappa+m\right)  p_{m-1,n}+\kappa
p_{m,n-1}-\kappa p_{n-1,m}\\
&  -\frac{\kappa\left(  m-n+1\right)  }{\kappa+m-n+1}\left(  N\kappa
+m+1\right)  p_{m,n-1}+g_{1},
\end{align*}
then since $\omega_{m-1,n}=p_{m-1,n}-\frac{\kappa}{\kappa+m-n}\left(  \left(
m-n\right)  p_{m,n-1}+p_{n-1,m}\right)  +g_{2}$ it follows that $\mathcal{D}%
_{1}\omega_{mn}-\left(  N\kappa+m\right)  \omega_{m-1,n}=\frac{\left(  \left(
N-1\right)  \kappa+n\right)  \kappa}{\kappa+m-n}\left(  p_{n-1,m}-\frac
{\kappa}{\kappa+m-n+1}p_{m,n-1}\right)  +g_{3}$. This proves part (1).

Next suppose $m\geq n$ then
\begin{align*}
\mathcal{D}_{2}\omega_{mn}  &  =\left(  \left(  N-1\right)  \kappa+n\right)
p_{m,n-1}+\kappa p_{n-1,m}\\
&  -\frac{\kappa}{\kappa+m-n+1}\left(  N\kappa+m+1\right)  p_{n-1,m}+g_{4}\\
&  =\left(  \left(  N-1\right)  \kappa+n\right)  \left(  p_{m,n-1}%
-\frac{\kappa}{\kappa+m-n+1}p_{n-1,m}\right)  +g_{4},
\end{align*}
and this proves part (2). Part (3) follows from (2) by setting $m=n$ and
applying the transposition $\left(  1,2\right)  .$
\end{proof}

The following evaluation formula for $\omega_{mn}\left(  1^{N}\right)  $
(where $1^{N}=\left(  1,1,\ldots,1\right)  \in\mathbb{R}^{N}$) is a special
case of a general result for nonsymmetric Jack polynomials (see \cite[p.310]%
{DX}). Here is a self-contained proof.

\begin{proposition}
\label{val1N}For $m\geq n,$%
\[
\omega_{mn}\left(  1^{N}\right)  =\frac{\left(  N\kappa+1\right)  _{m}\left(
\left(  N-1\right)  \kappa+1\right)  _{n}}{\left(  m-n\right)  !\,n!\,\left(
\kappa+m-n+1\right)  _{n}}.
\]
\end{proposition}

\begin{proof}
By the negative binomial theorem $p_{ij}\left(  1^{N}\right)  =\frac{\left(
N\kappa+1\right)  _{i}\left(  N\kappa+1\right)  _{j}}{i!\,j!}$. Substituting
this in Definition \ref{wdef} yields
\begin{gather*}
\omega_{mn}\left(  1^{N}\right)  =\frac{\left(  N\kappa+1\right)  _{m}\left(
N\kappa+1\right)  _{n}}{m!\,n!}\\
+\sum_{j=1}^{n}\frac{\left(  -\kappa\right)  _{j}\left(  m-n+1\right)  _{j-1}%
}{\left(  \kappa+m-n+1\right)  _{j}\,j!}\left(  m-n+2j\right)  \frac{\left(
N\kappa+1\right)  _{m+j}\left(  N\kappa+1\right)  _{n-j}}{\left(  m+j\right)
!\,\left(  n-j\right)  !}.
\end{gather*}
For now, assume $m>n$ then $\left(  m-n+1\right)  _{j-1}\left(  m-n+2j\right)
=$\newline $\left(  m-n\right)  _{j}\left(  \frac{m-n}{2}+1\right)
_{j}/\left(  \frac{m-n}{2}\right)  _{j}$ and
\begin{align*}
\omega_{mn}\left(  1^{N}\right)   &  =\frac{\left(  N\kappa+1\right)
_{m}\left(  N\kappa+1\right)  _{n}}{m!\,n!}\\
&  \times\sum_{j=0}^{n}\frac{\left(  -n\right)  _{j}\left(  -\kappa\right)
_{j}\left(  N\kappa+1+m\right)  _{j}\left(  m-n\right)  _{j}\left(  \frac
{m-n}{2}+1\right)  _{j}}{\left(  m+1\right)  _{j}\left(  \kappa+m-n+1\right)
_{j}\left(  -N\kappa-n\right)  _{j}\left(  \frac{m-n}{2}\right)  _{j}\,j!}.
\end{align*}
The sum is a terminating well-poised $_{5}F_{4}$ whose value is
\[
\frac{\left(  m-n+1\right)  _{n}\left(  -N\kappa+\kappa-n\right)  _{n}%
}{\left(  \kappa+m-n+1\right)  _{n}\left(  -N\kappa-n\right)  _{n}},
\]
a formula of Dougall (see Bailey \cite[p.25]{B}). The stated formula follows
by using the reversal $\left(  a-n\right)  _{n}=\left(  -1\right)  ^{n}\left(
1-a\right)  _{n}$. The formula is also valid when $m=n$ (consider $z=m-n$ as a
variable, then the limit of $\left(  z\right)  _{j}\left(  \frac{z}%
{2}+1\right)  _{j}/\left(  \frac{z}{2}\right)  _{j}$ as $z\rightarrow0$ is $1$
for $j=0$ and $2\times j!$ for $j\geq1$).
\end{proof}

\section{Restriction to $N=2$}

The main results of this paper revolve around the vanishing of $\omega_{mn}$
for certain values of $N,m,n$ and $\kappa$, but it is also necessary to show
certain $\omega_{mn}\neq0$. This will be accomplished by restricting to $N=2$
(setting $x_{i}=0$ for $i>2$) and finding an explicit formula for $\omega_{mn}.$

\begin{definition}
For $m\geq n$ let
\[
f_{mn}\left(  x_{1},x_{2}\right)  =\left(  x_{1}x_{2}\right)  ^{n}\sum
_{j=0}^{m-n}\frac{\left(  \kappa+1\right)  _{m-n-j}\left(  \kappa\right)
_{j}}{\left(  m-n-j\right)  !\,j!}x_{1}^{m-n-j}x_{2}^{j}.
\]
\end{definition}

\begin{proposition}
For $m\geq n$, $\mathcal{D}_{1}x_{1}f_{mn}=\left(  \kappa+m+1\right)  f_{mn}$
and \newline $\left(  \mathcal{D}_{2}x_{2}-\kappa\left(  1,2\right)  \right)
f_{mn}=\left(  n+1\right)  f_{mn}.$
\end{proposition}

\begin{proof}
It is clear from the generating function (\ref{genfun}) that $\omega
_{m,0}=p_{m,0}=f_{m,0}$ for $m\geq0$. This implies $\mathcal{D}_{1}%
x_{1}f_{m,0}=\left(  \kappa+m+1\right)  f_{m,0}$ and $\mathcal{D}_{2}%
x_{2}f_{m,0}=\left(  1+\kappa\left(  1,2\right)  \right)  f_{m,0}$ (it can be
shown directly from the generating function that $\mathcal{D}_{1}x_{1}%
p_{m,0}=\left(  \kappa+m+1\right)  p_{m,0}$ for $N=2$). By the product rule
(\ref{prodrule})
\begin{align*}
\mathcal{D}_{1}x_{1}f_{mn} &  =\left(  x_{1}x_{2}\right)  ^{n}\mathcal{D}%
_{1}x_{1}f_{m-n,0}+x_{1}f_{m-n,0}\frac{\partial}{\partial x_{1}}\left(
x_{1}x_{2}\right)  ^{n}\\
&  =\left(  \kappa+m-n+1+n\right)  \left(  x_{1}x_{2}\right)  ^{n}f_{m-n,0},
\end{align*}
and
\begin{align*}
\left(  \mathcal{D}_{2}x_{2}-\kappa\left(  1,2\right)  \right)  f_{mn} &
=\left(  x_{1}x_{2}\right)  ^{n}\left(  \mathcal{D}_{2}x_{2}-\kappa\left(
1,2\right)  \right)  f_{m-n,0}+x_{2}f_{m-n,0}\frac{\partial}{\partial x_{2}%
}\left(  x_{1}x_{2}\right)  ^{n}\\
&  =\left(  n+1\right)  \left(  x_{1}x_{2}\right)  ^{n}f_{m-n,0},
\end{align*}
as claimed.
\end{proof}

Since the joint eigenfunctions of the commuting operators $\mathcal{D}%
_{1}x_{1}$ and \newline $\left(  \mathcal{D}_{2}x_{2}-\kappa\left(
1,2\right)  \right)  $ are uniquely determined for generic $\kappa$ (including
$\kappa>0$) we see that $f_{mn}$ is a scalar multiple of $\omega_{mn}$.
Evaluation at $x=\left(  1,1\right)  $ determines the constant.

\begin{proposition}
For $N=2,m\geq n,$%
\[
\omega_{mn}=\frac{\left(  2\kappa+m-n+1\right)  _{n}\left(  \kappa+1\right)
_{n}}{\left(  \kappa+m-n+1\right)  _{n}\,n!}f_{mn}.
\]
\end{proposition}

\begin{proof}
By Proposition \ref{val1N} $\omega_{mn}\left(  1^{2}\right)  =\dfrac{\left(
2\kappa+1\right)  _{m}\left(  \kappa+1\right)  _{n}}{\left(  \kappa
+m-n+1\right)  _{n}\left(  m-n\right)  !n!}$ while $f_{mn}\left(
1^{2}\right)  =\dfrac{\left(  2\kappa+1\right)  _{m-n}}{\left(  m-n\right)
!}$ (by the Vandermonde sum formula). Since $\left(  2\kappa+1\right)
_{m}=\left(  2\kappa+1\right)  _{m-n}\left(  2\kappa+m-n+1\right)  _{n}$ this
completes the proof.
\end{proof}

We observe that $f_{mn}\neq0$ provided (as assumed throughout) that
$\kappa\notin-\mathbb{N}$. This leads to the following nontriviality result.

\begin{corollary}
\label{non0}For $N\geq2,m\geq n$ the polynomial $\omega_{mn}\neq0$ provided
that $2\kappa\neq-j$ where $j=m-n+1,m-n+2,\ldots,m$.
\end{corollary}

\section{Some singular polynomials}

These results already appeared in \cite{D2}, and serve as illustration.

\begin{proposition}
For $N\geq2$ and $n\in\mathbb{N}$ such that $N\nmid n$ (equivalently,
$\gcd\left(  N,n\right)  <N$), $\omega_{n,0}$ is a singular polynomial for
$\kappa=-\frac{n}{N}$.
\end{proposition}

\begin{proof}
Note that $\omega_{n,0}=p_{n,0}$. By formula (\ref{D1p}) $\mathcal{D}%
_{1}\omega_{n,0}=\left(  N\kappa+n\right)  \omega_{n-1,0}$ and $\mathcal{D}%
_{i}\omega_{n,0}=0$ for each $i>1$. Further $p_{n,0}\left(  1,0,\ldots\right)
=\frac{\left(  \kappa+1\right)  _{n}}{n!}$ which is not zero provided
$\kappa\notin-\mathbb{N}$.
\end{proof}

These polynomials have been studied by Chmutova and Etingof \cite{CE} in the
context of representations of the rational Cherednik algebra.

\begin{proposition}
For $N\geq4$ and $n\in\mathbb{N}$ such that $\gcd\left(  N-1,n\right)
<\frac{N-1}{2}$, $\omega_{nn}$ is a singular polynomial for $\kappa=-\frac
{n}{N-1}$.
\end{proposition}

\begin{proof}
The condition $\gcd\left(  N-1,n\right)  <\frac{N-1}{2}$ is equivalent to
excluding the values $\kappa=-j,-j+\frac{1}{2}$ for $j\in\mathbb{N}$. By
Theorem \ref{Dwmn} we have $\mathcal{D}_{1}\omega_{nn}=\left(  \left(
N-1\right)  \kappa+n\right)  \left(  \omega_{n-1,n}-\dfrac{\kappa}{\kappa
+1}\omega_{n,n-1}\right)  $ which is zero for $\kappa=-\frac{n}{N-1}$,
similarly $\mathcal{D}_{2}\omega_{nn}=0$, and $\mathcal{D}_{i}\omega_{nn}=0 $
for all $i>2$ (for all $\kappa).$ Corollary \ref{non0} shows that $\omega
_{nn}\neq0$ since $2\kappa\notin-\mathbb{N}$.
\end{proof}

Suppose that $M$ is an irreducible $S_{N}$-module of homogeneous polynomials
(that is, $M$ is a linear subspace of the space of polynomials, and is
invariant under the action of each $w\in S_{N}$, and has no proper nontrivial
invariant subspaces) then $M$ is of some isotype (corresponding to an
irreducible representation of $S_{N}$) labeled by a partition $\tau$ of $N$.

\begin{proposition}
\label{mutau}Suppose $\tau$ is a partition of $N$ and the homogeneous
polynomial $f$ is of degree $n$ and of isotype $\tau$, that is, $\mathrm{span}%
\left\{  wf:w\in S_{N}\right\}  $ is an $S_{N}$-module on which $S_{N}$ acts
by the irreducible representation corresponding to $\tau$, then $\sum
_{i=1}^{N}x_{i}\mathcal{D}_{i}f=\left(  n+\kappa\mu\left(  \tau\right)
\right)  f$, where $\mu\left(  \tau\right)  =\binom{N}{2}-\frac{1}{2}%
\sum_{j=1}^{N}\tau_{j}\left(  \tau_{j}+1-2j\right)  $.
\end{proposition}

\begin{proof}
It is easy to show that $\sum\limits_{i=1}^{N}x_{i}\mathcal{D}_{i}%
f=\sum\limits_{i=1}^{N}x_{i}\frac{\partial f}{\partial x_{i}}+\kappa
\sum\limits_{1\leq i<j\leq N}\left(  1-\left(  i,j\right)  \right)  f$ for any
polynomial $f$. The operator $\sum\limits_{1\leq i<j\leq N}\left(  1-\left(
i,j\right)  \right)  $ is constant on \newline $\mathrm{span}\left\{  wf:w\in
S_{N}\right\}  $, and its value is given by Young's formula,\newline
$\binom{N}{2}-\frac{1}{2}\sum_{j=1}^{N}\tau_{j}\left(  \tau_{j}+1-2j\right)  $
(see \cite[p.177]{D1}). The Euler operator $\sum\limits_{i=1}^{N}x_{i}%
\frac{\partial}{\partial x_{i}}$gives the degree of $f$.
\end{proof}

Thus a necessary condition for a homogeneous polynomial $f$ of isotype $\tau$
to be singular is that $\kappa=-\frac{\deg f}{\mu\left(  \tau\right)  }.$ The
isotype for $\mathrm{span}\left\{  w\omega_{n,0}:w\in S_{N}\right\}  $ when
$\kappa=-\frac{n}{N}$ is $\left(  N-1,1\right)  $ and $\mu\left(  \left(
N-1,1\right)  \right)  =N$. For the values $\kappa=-\frac{n}{N-1}$with
$2\kappa\notin-\mathbb{N}$ the isotype of $\mathrm{span}\left\{  w\omega
_{nn}:w\in S_{N}\right\}  $ is $\left(  N-2,2\right)  $ and $\mu\left(
\left(  N-2,2\right)  \right)  =2N-2$. It is exactly the filling of the gap at
$\gcd\left(  N-1,n\right)  =\frac{N-1}{2}$ for $N$ being odd that we consider
in the sequel. Before we leave this section we point out that the singular
polynomials described so far do not depend on $N$ (with the exception just
noted), that is, the singularity property holds for all $N$. This no longer
holds once we consider isotypes corresponding to partitions with more than two parts.

\section{Singular polynomials for half-integer parameter values}

In this section we show that the polynomials $\omega_{\left(  2l+1\right)
\left(  m+1\right)  ,\left(  2l+1\right)  m}$ are singular for $N=2m+1$ and
$\kappa=-l-\frac{1}{2}$, for $l\in\mathbb{N}_{0}$ and $m\in\mathbb{N}.$ By
Theorem \ref{Dwmn} we already know $\mathcal{D}_{i}\omega_{\left(
2l+1\right)  \left(  m+1\right)  ,\left(  2l+1\right)  m}=0$ for $i\geq2$.
Thus we need to show that $\omega_{\left(  2l+1\right)  \left(  m+1\right)
-1,\left(  2l+1\right)  m}=0$ for these choices of $N,\kappa$. This will be
done by introducing a new basis of polynomials related to the $\left\{
p_{mn}\right\}  $ basis by a linear relation involving Krawtchouk polynomials.

\begin{definition}
\label{qdef}The homogeneous polynomials $q_{mn}$ (for $m,n\in\mathbb{N}_{0}$)
are defined by
\[
\sum_{m,n=0}^{\infty}q_{mn}\left(  x\right)  u^{m}v^{n}=\sum_{i,j=0}^{\infty
}p_{ij}\left(  x\right)  \left(  u+v\right)  ^{i}\left(  u-v\right)  ^{j},
\]
the generating function converges for $\left|  u\right|  ,\left|  v\right|
<\left(  \max_{i}\left|  x_{i}\right|  \right)  /2$.
\end{definition}

We state the basic properties of the symmetric Krawtchouk polynomials (see
Szeg\"{o},\cite[p.36]{Sz}). They are orthogonal for the binomial distributions
with parameter $\frac{1}{2}$. Fix $n\in\mathbb{N}$ then the Krawtchouk
polynomial of degree $m$ (parameters $n,\frac{1}{2}$) with $0\leq m\leq n$ is
given by
\[
K_{m}\left(  t;n\right)  =\frac{1}{\binom{n}{m}}\sum_{j=0}^{m}\frac{\left(
t-n\right)  _{m-j}\left(  -t\right)  _{j}}{\left(  m-j\right)  !\,j!}\left(
-1\right)  ^{m-j}.
\]
Then the following hold for $0\leq m,l\leq n$:

\begin{enumerate}
\item $K_{m}\left(  0;n\right)  =1$, normalization;

\item $\left(  1-s\right)  ^{l}\left(  1+s\right)  ^{n-l}=\sum_{m=0}^{n}%
s^{m}\binom{n}{m}K_{m}\left(  l;n\right)  $, generating function;

\item $2^{-n}\sum_{t=0}^{n}\binom{n}{t}K_{m}\left(  t;n\right)  K_{l}\left(
t;n\right)  =\delta_{ml}\binom{n}{m}^{-1}$, orthogonality;

\item $K_{m}\left(  l;n\right)  =\,_{2}F_{1}\left(  -l,-m;-n;2\right)  $,
hypergeometric polynomial;

\item $K_{m}\left(  l;n\right)  =K_{l}\left(  m;n\right)  $, symmetry;

\item $K_{m}\left(  n-t;n\right)  =\left(  -1\right)  ^{m}K_{m}\left(
t;n\right)  $, parity$.$
\end{enumerate}

Any expansion in $\left\{  p_{n-i,i}:1\leq i\leq n\right\}  $ can be
transformed to one in \newline $\left\{  q_{n-i,i}:1\leq i\leq n\right\}  $
($n\in\mathbb{N}$) by means of Krawtchouk polynomials.

\begin{lemma}
Suppose $n\in\mathbb{N}$ and $f=\sum_{i=0}^{n}c_{i}p_{n-i,i}$ with
coefficients $c_{i}\in\mathbb{Q}\left(  \kappa\right)  $, then
\[
f=\frac{1}{2^{n}}\sum_{i=0}^{n}q_{n-i,i}\sum_{j=0}^{n}\binom{n}{j}c_{j}%
K_{i}\left(  j;n\right)  .
\]
\end{lemma}

\begin{proof}
In Definition \ref{qdef} replace $u,v$ by $\frac{s+t}{2},\frac{s-t}{2}$
respectively then $p_{n-j,j}$ equals the coefficient of $s^{n-j}t^{j}$ in
$2^{-n}\sum_{i=0}^{n}q_{n-i,i}\left(  s+t\right)  ^{n-i}\left(  s-t\right)
^{i}=$\newline $2^{-n}\sum_{i=0}^{n}q_{n-i,i}\sum_{j=0}^{n}\binom{n}{j}%
K_{j}\left(  i;n\right)  s^{n-j}t^{j}$. The lemma now follows from the
symmetry relation.
\end{proof}

We use the lemma to show that for special values of $n,i,\kappa$ the
coefficients of $\omega_{n-i,i}$ with respect to $\left\{  q_{n-j,j}:0\leq
j\leq n\right\}  $ have a vanishing property: $\omega_{n-i,i}=\sum_{j=0}%
^{n}c_{j}q_{n-j,j}$ and $c_{j}=0$ for $j>n-2i$. We will also find a similar
result for $\omega_{n-i,i}+\omega_{i,n-i}$.

\begin{proposition}
\label{w2q0}Suppose $n$ is even, $1\leq i\leq\frac{n}{2}$ and $\kappa
=-\frac{1}{2}\left(  n-2i+1\right)  $, then $\omega_{n-i,i}=\sum
\limits_{j=0}^{n-2i}c_{j}q_{n-j,j}$ with coefficients $c_{j}\in\mathbb{Q}$.
\end{proposition}

\begin{proof}
Substitute $\kappa=-\frac{1}{2}\left(  n-2i+1\right)  $ in Definition
\ref{wdef} for $\omega_{n-i,i}$ to obtain
\[
\omega_{n-i,i}=p_{n-i,i}+\sum_{l=1}^{i}\left(  \frac{\left(  n-2i+1\right)
_{l}}{l!}p_{n-i+l,i-l}+\frac{\left(  n-2i+1\right)  _{l-1}}{\left(
l-1\right)  !}p_{i-l,n-i+l}\right)  .
\]
Extract the coefficients of $\omega_{n-i,i}$ with respect to $p_{n-j,j}$ as
follows: for $0\leq j\leq i-1$ replace $l$ by $i-j$ in the first part of the
sum, the value is
\[
\frac{\left(  n-2i+1\right)  _{i-j}}{\left(  i-j\right)  !}=\frac{\left(
n-2i+i-j\right)  !}{\left(  n-2i\right)  !\left(  i-j\right)  !}=\frac{\left(
i-j+1\right)  _{n-2i}}{\left(  n-2i\right)  !};
\]
this is also valid for $j=i$; the coefficient is zero for $i+1\leq j\leq n-i$;
for $n-i<j\leq n$ replace $l$ by $j-n+i$ then
\[
\frac{\left(  n-2i+1\right)  _{l-1}}{\left(  l-1\right)  !}=\frac{\left(
l\right)  _{n-2i}}{\left(  n-2i\right)  !}=\frac{\left(  j-n+i\right)
_{n-2i}}{\left(  n-2i\right)  !}=\left(  -1\right)  ^{n}\frac{\left(
i-j+1\right)  _{n-2i}}{\left(  n-2i\right)  !}%
\]
(reversal of the Pochhammer symbol) but $n$ is even, so the value is
$\frac{\left(  i-j+1\right)  _{n-2i}}{\left(  n-2i\right)  !}$ just as for
$0\leq j\leq i.$ The expression $\left(  i-j+1\right)  _{n-2i}$ is a
polynomial in $j$, vanishing at $i+1,i+2,\ldots,n-i$. Thus
\begin{align*}
\omega_{n-i,i}  &  =\sum_{j=0}^{n}\frac{\left(  i-j+1\right)  _{n-2i}}{\left(
n-2i\right)  !}p_{n-j,j}\\
&  =\frac{1}{2^{n}}\sum_{l=0}^{n}q_{n-l,l}\sum_{j=0}^{n}\binom{n}{j}%
\frac{\left(  i-j+1\right)  _{n-2i}}{\left(  n-2i\right)  !}K_{l}\left(
j;n\right)  ,
\end{align*}
by the Lemma. The orthogonality property of $K_{l}$ shows that the coefficient
of $q_{n-l,l}$ vanishes for $l>n-2i$.
\end{proof}

\begin{proposition}
\label{w2q1}Suppose $n$ is odd, $0\leq i<\frac{n}{2}$ and $\kappa=-\frac{1}%
{2}\left(  n-2i\right)  $, then $\omega_{n-i,i}+\omega_{i,n-i}=\sum
_{j=0}^{n-2i-1}c_{j}q_{n-j,j}$ with coefficients $c_{j}\in\mathbb{Q}$.
\end{proposition}

\begin{proof}
It follows from Definition \ref{wdef} that
\[
\omega_{n-i,i}+\omega_{i,n-i}=\sum_{l=0}^{i}\frac{\left(  -\kappa\right)
_{l}\left(  n-2i\right)  _{l}}{\left(  \kappa+n-2i+1\right)  _{l}\,l!}%
\frac{n-2i+2l}{n-2i}\left(  p_{n-i+l,l-i}+p_{l-i,n-i+l}\right)  .
\]
When $\kappa=-\frac{1}{2}\left(  n-2i\right)  $ we obtain $\frac{\left(
-\kappa\right)  _{l}}{\left(  \kappa+n-2i+1\right)  _{l}\,}=\frac
{n-2i}{n-2i+2l}$. As before, $\frac{\left(  n-2i\right)  _{l}}{l!}%
=\frac{\left(  l+1\right)  _{n-2i-1}}{\left(  n-2i-1\right)  !}$, and replace
$l$ by $i-j$ and $j-n+i$ respectively for the ranges $0\leq j\leq i$ and
$n-i\leq i\leq n$ respectively. Also $\left(  j-n+i+1\right)  _{n-2i-1}%
=\left(  -1\right)  ^{n-1}\left(  i-j+1\right)  _{n-2i-1}$ and $n$ is odd. The
polynomial $\left(  i-j+1\right)  _{n-2i-1}$ in $j$ vanishes at $i+1,\ldots
,n-i-1$. Similarly to the previous proposition we find that
\[
\omega_{n-i,i}+\omega_{i,n-i}=\frac{1}{2^{n}}\sum_{l=0}^{n}q_{n-l,l}\sum
_{j=0}^{n}\binom{n}{j}\frac{\left(  i-j+1\right)  _{n-2i-1}}{\left(
n-2i-1\right)  !}K_{l}\left(  j;n\right)  ,
\]
and the coefficient of $q_{n-l,l}$ vanishes for $l>n-2i-1$. Because $\left(
i-j+1\right)  _{n-2i-1}=\left(  i-\left(  n-j\right)  +1\right)  _{n-2i-1}$
the coefficients of $q_{n-l,l}$ also vanish when $l$ is odd.
\end{proof}

We finish the construction of singular polynomials by showing for certain
values of $\kappa,N,n,i$ that the polynomials $q_{n-l,l}$ vanish for $0\leq
l\leq n-2i$. First we consider some partial products in the generating function.

\begin{definition}
The power series $A_{n}\left(  u;\kappa\right)  ,B_{n}\left(  u;\kappa\right)
$ (arbitrary $\kappa\in\mathbb{C},\,\left|  u\right|  ,\left|  v\right|
<\frac{1}{2}$ and $n\in\mathbb{N}_{0}$) are given by
\begin{align*}
\left(  1-\left(  u+v\right)  \right)  ^{-\kappa}\left(  1-\left(  u-v\right)
\right)  ^{-\kappa}  &  =\sum_{n=0}^{\infty}A_{n}\left(  u;\kappa\right)
v^{n},\\
\left(  1-\left(  u+v\right)  \right)  ^{-\kappa-1}\left(  1-\left(
u-v\right)  \right)  ^{-\kappa}  &  =\sum_{n=0}^{\infty}B_{n}\left(
u;\kappa\right)  v^{n}.
\end{align*}
\end{definition}

When $\kappa=-l-\frac{1}{2},l\in\mathbb{N}_{0}$ some of the series
$A_{n},B_{n}$ are actually polynomials.

\begin{lemma}
For $n\leq2l$, the functions $A_{n}\left(  u;-l-\frac{1}{2}\right)
,B_{n}\left(  u;-l-\frac{1}{2}\right)  $ are polynomials in $u$ of degree
$2l+1-n,2l-n$ respectively.
\end{lemma}

\begin{proof}
For the first part
\begin{align*}
\sum_{n=0}^{\infty}A_{n}\left(  u;\kappa\right)  v^{n}  &  =\left(
1-u\right)  ^{-2\kappa}\left(  1-\left(  \frac{v}{1-u}\right)  ^{2}\right)
^{-\kappa}\\
&  =\sum_{j=0}^{\infty}\frac{\left(  \kappa\right)  _{j}}{j!}v^{2j}\left(
1-u\right)  ^{-2\kappa-2j}.
\end{align*}
Thus $A_{n}=0$ if $n$ is odd and $A_{2j}\left(  u;-l-\frac{1}{2}\right)
=\frac{\left(  -l-\frac{1}{2}\right)  _{j}}{j!}\left(  1-u\right)  ^{2l+1-2j}%
$, which is a polynomial of degree $2l+1-2j$ provided $2j\leq2l$. For the
second part
\begin{gather*}
\sum_{n=0}^{\infty}B_{n}\left(  u;\kappa\right)  v^{n}=\left(  1-u\right)
^{-2\kappa-1}\left(  1+\frac{v}{1-u}\right)  \left(  1-\left(  \frac{v}%
{1-u}\right)  ^{2}\right)  ^{-\kappa-1}\\
=\sum_{j=0}^{\infty}\frac{\left(  \kappa+1\right)  _{j}}{j!}\left(
v^{2j}\left(  1-u\right)  ^{-2\kappa-1-2j}+v^{2j+1}\left(  1-u\right)
^{-2\kappa-2j-2}\right)  .
\end{gather*}
Thus $B_{n}\left(  u;-l-\frac{1}{2}\right)  =\frac{\left(  -l+\frac{1}%
{2}\right)  _{j}}{j!}\left(  1-u\right)  ^{2l-n}$ (with $j=\left\lfloor
\frac{n}{2}\right\rfloor $), a polynomial of degree $2l-n$ provided $n\leq2l$.
\end{proof}

To be precise the polynomials $A_{n}\left(  u;-l-\frac{1}{2}\right)  $ are of
degree $\leq2l+1-n$ (being $0$ when $n$ is odd).

\begin{proposition}
\label{q2z}Let $\kappa=-l-\frac{1}{2},l\in\mathbb{N}_{0}$ and suppose that
$n\leq2l$ then $m+n\geq N\left(  2l+1\right)  -1$ implies $q_{mn}=0$.
\end{proposition}

\begin{proof}
By combining the generating function (\ref{genfun}) for $\left\{
p_{mn}\right\}  $ and Definition \ref{qdef} we obtain
\begin{align*}
\sum_{m,n=0}^{\infty}q_{mn}u^{m}v^{n} &  =\sum_{\alpha\in\mathbb{N}_{0}^{N}%
}B_{\alpha_{1}}\left(  ux_{1};\kappa\right)  \left(  -1\right)  ^{\alpha_{2}%
}B_{\alpha_{2}}\left(  ux_{2};\kappa\right)  \\
&  \times\prod_{s=3}^{N}A_{\alpha_{s}}\left(  ux_{s};\kappa\right)  x^{\alpha
}v^{\left|  \alpha\right|  }.
\end{align*}
For a fixed $n\leq2l$ the coefficient of $v^{n}$ is the sum over $\alpha
\in\mathbb{N}_{0}^{N}$ with $\left|  \alpha\right|  =n$. For each $\alpha$
with $\left|  \alpha\right|  =n$ (implying that each $\alpha_{i}\leq2l$ and
the Lemma applies) the corresponding term is a product of polynomials in $u$
of degree $\left(  2l-\alpha_{1}\right)  +\left(  2l-\alpha_{2}\right)
+\sum_{s=3}^{N}\left(  2l+1-\alpha_{s}\right)  =N\left(  2l+1\right)  -2-n$.
This shows that the coefficient of $u^{m}v^{n}$ vanishes if $m\geq N\left(
2l+1\right)  -1-n$.
\end{proof}

\begin{theorem}
Let $\kappa=-l-\frac{1}{2}$ and $N=2m+1$, with $l\in\mathbb{N}_{0}%
,m\in\mathbb{N}$ then $\omega_{\left(  2l+1\right)  \left(  m+1\right)
,\left(  2l+1\right)  m}$ is a singular polynomial, of isotype $\left(
N-2,1,1\right)  $, and \newline $\left(  1+\left(  1,2\right)  \right)
\omega_{\left(  2l+1\right)  \left(  m+1\right)  ,\left(  2l+1\right)  m}=0$.
\end{theorem}

\begin{proof}
Let $a=\left(  2l+1\right)  \left(  m+1\right)  ,b=\left(  2l+1\right)  m$.
Since $-2\kappa=2l+1<a-b+1$ Corollary \ref{non0} shows that $\omega_{ab}\neq
0$. By Theorem \ref{Dwmn} $\mathcal{D}_{2}\omega_{ab}=\left(  2m\kappa
+b\right)  \left(  \omega_{a,b-1}-\frac{\kappa}{\kappa+a-b+1}\omega
_{b-1,a}\right)  =0 $ for $\kappa=-\frac{b}{2m}=-l-\frac{1}{2}.$ Similarly
$\mathcal{D}_{1}\omega_{ab}=\left(  N\kappa+a\right)  \omega_{a-1,b}.$ By
Proposition \ref{w2q0} $\omega_{a-1,b}=\sum_{j=0}^{a-1-b}c_{j}q_{a+b-1-j,j}$
(since -$\frac{1}{2}\left(  a-b\right)  =\kappa$) with some coefficients
$c_{j}\in\mathbb{Q}$. Since $a-1-b=2l$ and $a+b-1=N\left(  2l+1\right)  -1$,
Proposition \ref{q2z} shows that each $q_{a+b-1-j,j}=0$ (for $j\leq2l$).

Let $M=\mathrm{span}\left\{  w\omega_{\left(  2l+1\right)  \left(  m+1\right)
,\left(  2l+1\right)  m}:w\in S_{N}\right\}  $. By the invariance properties
of $\left\{  \mathcal{D}_{i}:1\leq i\leq N\right\}  $ any nonzero element of
$M$ is also singular. In general for $m>n$ the $S_{N}$-module $\mathrm{span}%
\left\{  w\omega_{mn}:w\in S_{N}\right\}  $, which realizes the representation
of $S_{N}$ induced up from the trivial representation of $S_{1}\times
S_{1}\times S_{N-2}$, decomposes into the isotypes $\left(  N-2,1,1\right)
,\left(  N-2,2\right)  ,\left(  N-1,1\right)  ,\left(  N\right)  $ (see
\cite[p.115]{M}). The eigenvalues $\mu\left(  \tau\right)  $ are $2N,2N-2,N,0$
respectively (see Proposition \ref{mutau}); but the singularity condition
implies $N\left(  2l+1\right)  +\kappa\mu\left(  \tau\right)  =0$ thus the
latter three can not contain singular polynomials for $\kappa=-l-\frac{1}{2}$
(note the degree of the polynomial $\omega_{ab}$ is $N\left(  2l+1\right)  $).
This implies that $M$ is of isotype $\left(  N-2,1,1\right)  $ and hence is of
dimension $\binom{N-1}{2}$.

By Proposition \ref{w2q1} $\omega_{ab}+\omega_{ba}=\sum_{j=0}^{a-b-1}%
c_{j}q_{a+b-j,j}$ (since $-\frac{1}{2}\left(  a-b\right)  =\kappa$ and $a+b$
is odd) with $c_{j}\in\mathbb{Q}$. Since $a-b-1=2l$ and $a+b=N\left(
2l+1\right)  $, Proposition \ref{q2z} shows that each $q_{a+b-j,j}=0$.
\end{proof}

It is interesting that the parameters of the singular polynomials just barely
satisfy the various inequalities appearing in the preparatory results.

\bigskip

\textsc{Department of Mathematics, University of Virginia}

\textsc{Charlottesville, VA 22904-4137, U.S.}

\textbf{cfd5z@virginia.edu}
\end{document}